\documentclass[11pt, a4paper]{amsart}
\usepackage{amsmath, amssymb, amsthm, amscd}
\usepackage{graphicx}

\theoremstyle{plain}
\newtheorem{theorem}{Theorem}
\newtheorem{lemma}[theorem]{Lemma}
\newtheorem{proposition}[theorem]{Proposition}
\newtheorem{corollary}[theorem]{Corollary}

\theoremstyle{definition}
\newtheorem{definition}{Definition}
\newtheorem{conjecture}{Conjecture}

\newcommand{\tr}{\operatorname{tr\,}}
\newcommand{\cl}{\operatorname{cl\,}}
\renewcommand{\Im}{\operatorname{Im\,}}
\renewcommand{\Re}{\operatorname{Re\,}}
\renewcommand{\int}{\operatorname{int}}
\newcommand{\Inv}{\operatorname{Inv}}
\newcommand{\Aut}{\operatorname{Aut}}
\newcommand{\Fix}{\operatorname{Fix}}

\title[On Loops in the Hyperbolic Locus]
{On Loops in the Hyperbolic Locus \\ 
 of the Complex H\'{e}non Map \\ 
 and Their Monodromies}

\author[Z. Arai]{Zin ARAI}

\address{Department of Mathematics, Kyoto University, Kyoto 606-8502, Japan}
\email{arai@math.kyoto-u.ac.jp}
\urladdr{http://www.math.kyoto-u.ac.jp/~arai/}
\thanks{Research supported in part by 
Grant-in-Aid for Scientific Research (No. 17740054), 
Japan Society for Promotion of Science.}

\begin{document}

\begin{abstract}
 We prove John Hubbard's conjecture on the topological complexity of
 the hyperbolic horseshoe locus of the complex H\'{e}non map. 
 Indeed, we show that there exist several non-trivial loops in the locus 
 which generate infinitely many  mutually different monodromies.
 Our main tool is a rigorous computational algorithm for verifying 
 the uniform hyperbolicity of chain recurrent sets.
 In addition, we show that the dynamics of the real H\'{e}non map is
 completely determined by the monodromy of a certain loop, providing 
 the parameter of the map is contained in the hyperbolic horseshoe 
 locus of the complex H\'{e}non map. 
\end{abstract}

\maketitle

\section{Introduction}

One of the motivations of this work is to give an answer to the conjecture
of John Hubbard on the topology of hyperbolic horseshoe locus of the complex
H\'{e}non map
\begin{equation*}
 H_{a,c}: \mathbb{C}^2 \to \mathbb{C}^2: 
   \begin{pmatrix} x \\ y \end{pmatrix} 
   \mapsto 
    \begin{pmatrix} x^2 +c - ay \\ x\end{pmatrix}.
\end{equation*}
Here $a$ and $c$ are complex parameters. 

We describe the conjecture following a formulation given by Bedford and 
Smillie \cite{BS05}.

Let us define
\begin{equation*}
 K_{a,c}^{\mathbb{C}}:= 
  \{p \in \mathbb{C}^2: \{H_{a,c}^n(p)\}_{n \in \mathbb{Z}} 
  \text{ is bounded}\},
  \quad K_{a,c}^{\mathbb{R}}:= K_{a,c}^{\mathbb{C}} \cap \mathbb{R}^2.
\end{equation*}
The set $K_{a, c}^{\mathbb{C}}$ is compact and invariant with respect to 
$H_{a,c}$.
When the parameters $a$ and $c$ are both real, 
the real plane $\mathbb{R}^2 \subset \mathbb{C}^2$
is invariant and hence so is $K_{a, c}^{\mathbb{R}}$.
In this case, we regard $H_{a, c}$ also as a dynamical system defined
on $\mathbb{R}^2$ and call it the real H\'{e}non map.

Our primary interest is on the parameter space, especially the set of 
parameters such that the complex and real H\'{e}non maps become a
uniformly hyperbolic horseshoe.
More precisely, we study the following sets:
\begin{gather*}
 \mathcal{H}^{\mathbb{C}} := \{(a,c) \in \mathbb{C}^2 : H_{a,c}|
 {K_{a,c}^{\mathbb{C}}} \text{ is a hyperbolic full horseshoe} \},\\
 \mathcal{H}^{\mathbb{R}} := \{(a,c) \in \mathbb{R}^2 : H_{a,c}|
 {K_{a,c}^{\mathbb{R}}} \text{ is a hyperbolic full horseshoe} \}.
\end{gather*}
By a hyperbolic full horseshoe, we mean an uniformly hyperbolic
invariant set which is topologically conjugate to the full shift map
$\sigma$ defined on $\Sigma_2 = \{0, 1\}^{\mathbb{Z}}$, the space of
bi-infinite sequences of two symbols.

A classical result of Devaney and Nitecki \cite{DN79}
claims that if $(a,c)$ is in
\begin{equation*}
 \mathrm{DN} := 
 \{(a,c) \in \mathbb{R}^2 : c < - (5+2\sqrt{5}) \, (|a|+1)^2/ 4, \, a \ne 0\}
\end{equation*}
then $K_{a, c}^{\mathbb{R}}$ is a hyperbolic full horseshoe.
Thus $\mathrm{DN} \subset \mathcal{H}^{\mathbb{R}}$ holds. 
They also showed that the set
\begin{equation*}
 \mathrm{EMP} := \{(a,c) \in \mathbb{R}^2 : c > (|a|+1)^2/ 4 \}
\end{equation*}
consists of parameter values such that $K_{a,c}^{\mathbb{R}} = \emptyset$.
Later, Hubbard and Oberste-Vorth investigated the H\'{e}non map
form the point of view of complex dynamics
and improved the hyperbolicity criterion
by showing that 
\begin{equation*}
  \mathrm{HOV} := \{(a,c) \in \mathbb{C}^2 : |c| > 2(|a|+1)^2, \, a \ne 0\}
\end{equation*}
is included in $\mathcal{H}^{\mathbb{C}}$.
Remark that $\mathrm{EMP} \cap \mathrm{HOV}$ is non-empty. 
In this parameter region, although $K_{a, c}^{\mathbb{C}}$ is a full 
horseshoe, it does not intersect with $\mathbb{R}^2$.

\begin{figure}[htb]
 \label{FIG:plateau}
 \begin{center}
  {\includegraphics{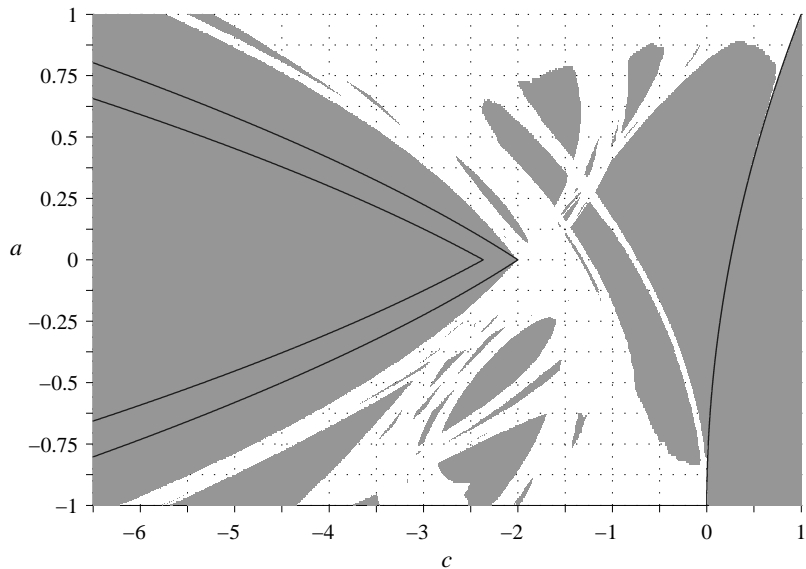}}
 \end{center}
 \vspace{-5mm}
 \caption{The shaded regions consist of
 hyperbolic (not necessarily full horseshoe) 
 parameters of the {\emph{real}} H\'{e}non map.}
\end{figure}

Figure \ref{FIG:plateau} illustrates a subset of parameter values 
on which the chain recurrent set of the {\it real} H\'{e}non map is uniformly
hyperbolic (not necessarily a full horseshoe) \cite{Arai05}.
Three solid lines are parts of the boundaries of $\mathrm{DN}$, $\mathrm{HOV}$
and $\mathrm{EMP}$, from left to right.
On the biggest island to the left, the chain recurrent set coincides
with $K_{a, c}^{\mathbb{R}}$ and is conjugate to the full shift.
Hence the island is contained in $\mathcal{H}^{\mathbb{R}}$.

We then consider the relation between $\mathcal{H}^{\mathbb{R}}$ and 
$\mathcal{H}^{\mathbb{C}}$.
By the result of Bedford, Lyubich and Smillie \cite[Theorem~10.1]{BLS}, 
we have
$\mathcal{H}^{\mathbb{R}} \subset \mathcal{H}^{\mathbb{C}} \cap \mathbb{R}^2$.
It is then natural to ask what happens in the rest of 
$\mathcal{H}^{\mathbb{C}} \cap \mathbb{R}^2$.

To be specific, we divide $\mathcal{H}^{\mathbb{C}} \cap \mathbb{R}^2$
into three mutually disjoint sets.

\begin{definition}[Bedford and Smillie \cite{BS05}]
 We call $(a,c) \in \mathcal{H}^{\mathbb{C}} \cap \mathbb{R}^2$ is 
 of type-1 if $(a,c) \in \mathcal{H}^{\mathbb{R}}$, and
 of type-2 if $K_{a,c}^{\mathbb{R}} = \emptyset$.
 Otherwise, it is of type-3.
\end{definition}

Since $\mathrm{DN} \subset \mathcal{H}^{\mathbb{R}}$, the set of type-1
parameter values is non-empty. 
The set of type-2 parameter values is also non-empty since it contains
$\mathrm{EMP} \cap \mathrm{HOV}$.
However, the existence of a type-3 parameter value was open.

\begin{conjecture}[Hubbard]
 There exists a parameter value of type-3.
\end{conjecture}

As we will see later, this conjecture turned to be true.

Besides the existence, Hubbard also conjectured
that there are infinitely many classes of type-3 parameter values
corresponding to mutually different real dynamics.
This stronger conjecture is, to be precise, given in terms of
the monodromy representation of the fundamental group of the hyperbolic
horseshoe locus as follows.

Denote by $\mathcal{H}^{\mathbb{C}}_0$ the component of 
$\mathcal{H}^{\mathbb{C}}$ that contains $\mathrm{HOV}$.
Let us fix a basepoint $(a_0, c_0) \in \mathrm{DN}$
and a topological conjugacy $h_0: K_{a_0,c_0}^{\mathbb{C}} \to \Sigma_2$.

Given a loop $\gamma: [0,1] \to \mathcal{H}^{\mathbb{C}}_0$
based at $(a_0,c_0)$, we construct a continuous family of
conjugacies $h_t^{\gamma}: K_{\gamma(t)}^{\mathbb{C}} \to \Sigma_2$
along $\gamma$ such that $h_0^{\gamma} = h_0$
(see \S \ref{SEC:monodromy} for the details).
This is possible because $K_{a, c}^{\mathbb{C}}$ is uniformly hyperbolic
along $\gamma$.
When no confusion may result we suppress $\gamma$ and write 
$h_t^{\gamma}$ as $h_t$.
Finally we set $\rho(\gamma) := h_1 \circ (h_0)^{-1}$.
It is easy to see that $\rho$ defines a group homomorphism
\begin{equation*}
  \rho: \pi_1 (\mathcal{H}^{\mathbb{C}}_0, (a_0, c_0)) \to \Aut(\Sigma_2)
\end{equation*}
where $\Aut(\Sigma_2)$ is the group of the automorphisms of 
$\Sigma_2$. 
Recall that an $automorphism$ of $\Sigma_2$ is a homeomorphism
of $\Sigma_2$ which commutes with the shift \cite{Kitchens}.
We call $\rho$ {\it the monodromy homomorphism} and denote its image by 
$\Gamma$.

For example, let $\gamma_{\emptyset}$ be a loop in 
$\mathcal{H}^{\mathbb{C}}_0$ based at $(a_0, c_0)$ and homotopic 
to the generator of $\pi_1(\mathrm{HOV})$.
It is then shown \cite{BS05} that $\rho(\gamma_{\emptyset})$ is 
an involution which interchanges the symbols $0$ and $1$.
Namely, $(\rho(\gamma_{\emptyset})(s))_i = 1 - s_i$ for 
all $s = (s_i) \in \{0, 1\}^{\mathbb{Z}}$.

The monodromy homomorphism was originally defined for polynomial
maps of one complex variable. 
In this case, since the map does not have the inverse, 
the target space of the monodromy homomorphism is the
automorphism group of {\it one-sided} shift space of $d$-symbols,
where $d$ is the degree of the polynomial. 
When $d = 2$, this group is isomorphic to $\mathbb{Z}_2$ and the monodromy
homomorphism is shown to be surjective since it maps the generator of 
$\pi_1(\mathbb{C} \setminus \{\text{the Mandelbrot set}\})$
to $1 \in \mathbb{Z}_2$.
The monodromy homomorphism is also surjective even when $d > 2$, although
the proof is much harder than the case $d = 2$ because the
automorphism group becomes much more complicated \cite{BDK1991}.

Hubbard conjectured that the surjectivity also holds in the case of 
the complex H\'{e}non map.

\begin{conjecture}[Hubbard]
 The monodromy homomorphism $\rho$ is surjective, that is, 
 $\Gamma = \Aut(\Sigma_2)$.
\end{conjecture}

The structure of $\Aut(\Sigma_2)$ is quite complicated \cite{BLR1988}:
it contains every finite group;
furthermore, it contains the direct sum of any countable collection of 
finite groups;
and it also contains the direct sum of countably many copies of $\mathbb{Z}$.
Therefore, the conjecture implies, provided it is true,
that the topological structure of $\mathcal{H}^{\mathbb{C}}$ is very rich, 
in contrast to the one-dimensional case where 
the fundamental group of 
$\mathbb{C} \setminus \{ \text{the Mandelbrot set}\}$ is simply $\mathbb{Z}$.

Let us state the main results of the paper now. 

First, we claim that Conjecture~1 is true.

\begin{theorem}
 \label{THM:existence}
 There exist parameter values of type-3.
 In fact, if $(a, c)$ is in one of the following sets:
 \begin{align*}
  &I_p := \{1\} \times [-5.46875, -5.3125],
  &I_q := \{0.25\} \times [-2.296875, -2.21875],\\
  &I_r := \{-1\} \times [-5.671875, -4.4375],
  &I_s := \{-0.375\} \times [-2.15625, -1.8125]
 \end{align*}
 then $(a, c)$ is of type 3.
\end{theorem}

As far as Conjecture~2 is concerned, we obtain the following result.

\begin{theorem}
 \label{THM:infinite_order}
 The order of the group $\Gamma$ is infinite.
 In particular, it contains an element of infinite order.
\end{theorem}

Apart form the theoretical interest, the monodromy theory 
of complex H\'{e}non map can contribute to the understanding of 
the real H\'{e}non map.

Let $(a, c) \in \mathcal{H}^{\mathbb{C}} \cap \mathbb{R}^2$.
If $(a, c)$
is of type-1 or 2, then by definition $K^{\mathbb{R}}_{a, c}$ is 
a full horseshoe, or the empty set.
Suppose $(a, c)$ is of type-3.
We then ask what $K^{\mathbb{R}}_{a, c}$ can be.
By definition, $K^{\mathbb{R}}_{a, c}$ is
a proper subset of $K_{a,c}^{\mathbb{C}} \cong \Sigma_2$.
The uniform hyperbolicity implies the existence of a Markov partition
for $K^{\mathbb{R}}_{a, c}$, and therefore, 
$K^{\mathbb{R}}_{a, c}$ must be topologically 
conjugate to some subshift of finite type.
The following theorem reveals that $K^{\mathbb{R}}_{a, c}$ is actually
a subshift of $\Sigma_2$ which is realized as the fixed point set of
the monodromy of a loop passing through $(a, c)$.

\begin{theorem}
 \label{THM:pruning}
 For any $(a, c) \in \mathcal{H}_0^{\mathbb{C}} \cap \mathbb{R}^2$,
 there exists a loop $\gamma: [0, 1] \to \mathcal{H}_0^{\mathbb{C}}$ with
 $\gamma(1/2) = (a, c)$ such that
 $H_{a, c}: K_{a,c}^{\mathbb{R}} \to K_{a,c}^{\mathbb{R}}$ is 
 topologically conjugate to
 \[
 \sigma|_{\Fix (\rho(\gamma))}: \Fix (\rho(\gamma)) \to \Fix (\rho(\gamma)).
 \]
 In fact, it suffice to set
 $\gamma := \bar{\alpha}^{-1} \cdot \alpha$,  where $\alpha$
 is an arbitrary path in $\mathcal{H}_0^{\mathbb{C}}$ 
 that starts at a point in $\mathrm{DN}$ and ends at $(a, c)$.
 Here $\bar{\alpha}$ denotes the complex conjugate of $\alpha$.
 The conjugacy is given by the restriction of $h_{1/2}$
 to $K^{\mathbb{R}}_{a, c}$. Namely, the following diagram commutes.
 \begin{equation*}
  \begin{CD}
   K^{\mathbb{R}}_{a, c}  @>{H_{a, c}}>> K^{\mathbb{R}}_{a, c}\\
   @V{h_{1/2}}V{\cong}V @V{\cong}V{h_{1/2}}V\\
   \Fix(\rho(\gamma)) @>\sigma>> \Fix(\rho(\gamma)).
  \end{CD}
 \end{equation*} 
\end{theorem}

As an application of Theorem~\ref{THM:pruning}, we obtain the following.

\begin{theorem}
 \label{THM:SFT}
 Let $(a, c) \in I_p$.
 The real H\'{e}non map 
 $H_{a, c}: K_{a,c}^\mathbb{R} \to  K_{a,c}^\mathbb{R}$ is topologically
 conjugate to the subshift of $\Sigma_2$ with two forbidden blocks
 $0010100$ and $0011100$.
 Similarly, $K_{a,c}^\mathbb{R}$ is conjugate to the subshift of $\Sigma_2$
 defined the following forbidden blocks: $10100$ and $11100$ for $(a, c) \in I_q$;
 $10010$ and $10110$ for $(a, c) \in I_r$;
 $0010$ and $0110$ for $(a, c) \in I_s$.
\end{theorem}

Notice that $I_p$ contains $(a, c) = (1, -5.4)$, the parameter studied by
Davis, MacKay and Sannami \cite{DMS91}. 
The subshift for $(a, c) \in I_p$
given in Theorem~\ref{THM:SFT} is equivalent to that 
observed by them.
Thus, we can say that their observation is now
rigorously verified.
We also remark that this theorem is closely related to the
so-called ``pruning front'' theory \cite{Cvi91, DMS91}.
Theorem~\ref{THM:pruning} implies that 
``primary pruned regions'', or, ``missing blocks'' of 
$K_{a,c}^{\mathbb{R}}$ is nothing else but
the region where the exchange of symbols occurs along $\gamma$.

The structure of the paper is as follows.
We prove the theorems in Section~\ref{SEC:proofs},
leaving computational algorithms to 
Section~\ref{SEC:hyperbolicity} and \ref{SEC:monodromy}.
In Section~\ref{SEC:hyperbolicity}, we summarize the algorithm 
for proving uniform hyperbolicity developed by the author \cite{Arai05}.
Section~\ref{SEC:monodromy} is devoted to an algorithm for computing 
the monodromy homomorphism.
In the appendix, we discuss a method for rigorously counting the number of
periodic points, which gives rise to an alternative proof of 
Theorem~\ref{THM:existence}.
Programs for computer assisted proofs are available at the author's 
pweb page (http://www.math.kyoto-u.ac.jp/\textasciitilde arai/).

The author is grateful, first of all, to John Hubbard, 
the originator of the problem.
He also would like to thank
E.~Bedford, P.~Cvitanovi\'{c}, S.~Hruska, H. Kokubu,
A.~Sannami, J.~Smillie, and S.~Ushiki
for many valuable suggestions.

\section{Proofs}\label{SEC:proofs}
We first prove Theorem~\ref{THM:pruning}.
The key is the symmetry of the H\'{e}non map with respect 
to the complex conjugation \cite{BS05}.
By the symmetry we mean the equation
\begin{equation*}
  \phi \circ H_{a, c} = H_{\bar{a}, \bar{c}} \circ \phi
\end{equation*}
where $\phi$ is the complex conjugation that
maps $z = (x, y)$ to $\bar{z} = (\bar{x}, \bar{y})$.

\begin{proof}[Proof of Theorem \ref{THM:pruning}]
 Take an arbitrary point $z \in K^{\mathbb{C}}_{a, c}$ and let
 \[
 s_z := h_{1/2}(z) \in \Sigma_2.
 \] 
 We will show that $\rho(\gamma) (s_z) = s_z$ if and only if 
 $z \in \mathbb{R}^2$.
 By abuse of notation, we denote the continuation of $z$ along $\gamma$ by
 \begin{equation*}
  z(\gamma, t) := (h_t)^{-1} (s_z) \in \mathbb{C}^2.
 \end{equation*} 
 Note that $z(\gamma, 1/2) = z$.
 By the continuity of hyperbolic invariant sets, the map 
 $t \mapsto (\gamma(t), z(\gamma, t))$
 defines a continuous loop in $\mathbb{C}^2 \times \mathbb{C}^2$.

 Let $\bar{\gamma} := \phi \circ \gamma$.
 The symmetry $\phi \circ H_{a, c} = H_{\bar{a}, \bar{c}} \circ \phi$ implies
 that the complex conjugate of the continuation of $z$ along $\gamma$ is just
 the continuation of $\bar{z} = \phi(z)$ along $\bar{\gamma}$.
 That is, we have
 \begin{equation*}
  \overline{z(\gamma, t)} = \bar{z}(\bar{\gamma},t).
 \end{equation*} 
 By the construction, we have $\bar{\gamma} = \gamma^{-1}$
 and hence $z(\bar{\gamma}, t) = z(\gamma, 1 -t)$.
 Therefore, 
 \begin{align*}
  \rho(\gamma) (s_z) 
  = h_1 ((h_0)^{-1} (s_z))
  = h_1 (z(\gamma, 0)) 
  = h_1 (\overline{z(\gamma, 0)})& \\
  = h_1 (\bar{z}(\bar{\gamma}, 0)) 
  = h_1 (\bar{z}(\gamma, 1)) 
  = h_1 ((h_1)^{-1}(s_{\bar{z}}))& = s_{\bar{z}}.
 \end{align*}
 The third equality holds because 
 $K_{\gamma(0)}^{\mathbb{C}} \subset \mathbb{R}^2$
 and hence $z(\gamma, 0) = \overline{z(\gamma, 0)}$.
 Since the map $h_{1/2}$ is a bijection between
 $K_{a,c}^{\mathbb{C}}$ and $\Sigma_2$, it follows that 
 $\rho(\gamma) (s_z) = s_z$ if and only if $z = \bar{z}$.
 This proves the theorem.
\end{proof}

\begin{figure}[htb]
 \begin{center}
  {\includegraphics{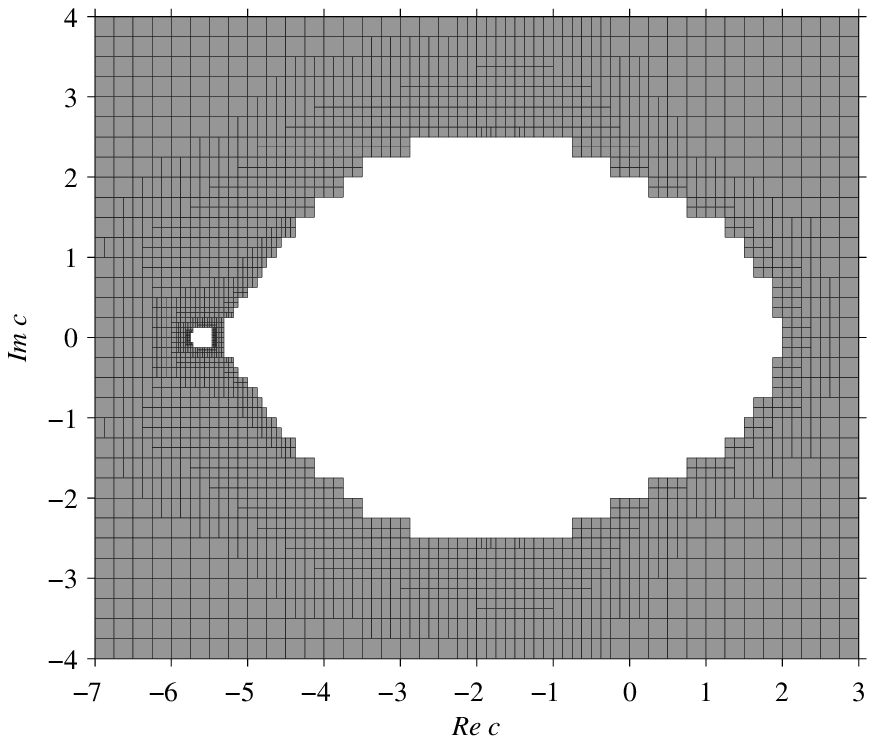}}
 \end{center}
 \caption{The shaded region is contained in 
 $\mathcal{H}^{\mathbb{C}} \cap \{a = 1\}$.}
 \label{FIG:ap}
\end{figure}

\begin{figure}[htb]
 \begin{center}
  {\includegraphics{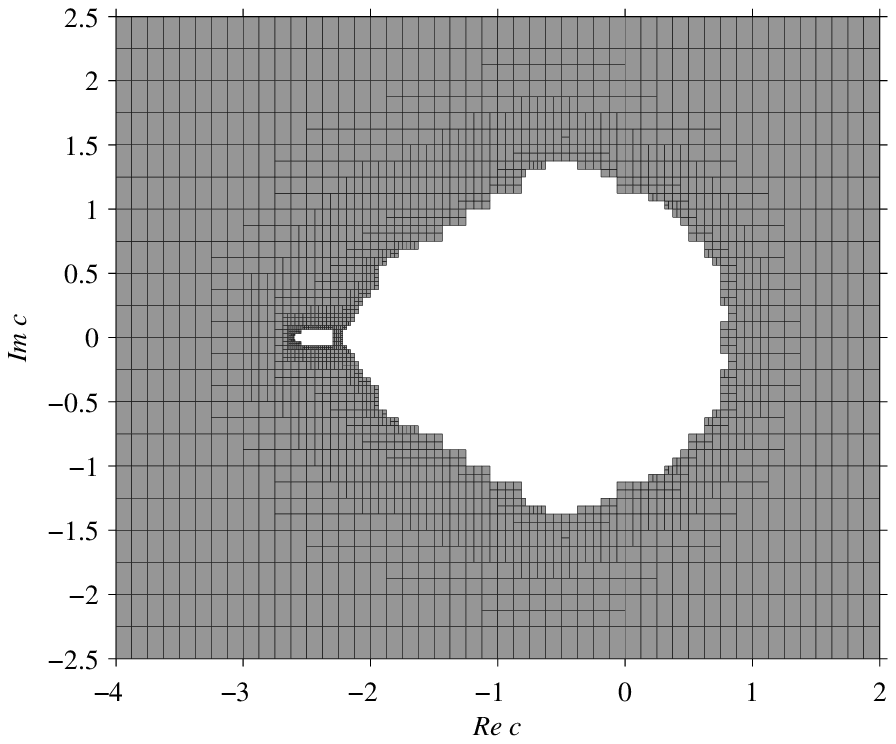}}
 \end{center}
 \caption{The shaded region is contained in 
 $\mathcal{H}^{\mathbb{C}} \cap \{a = 0.25\}$.}
 \label{FIG:a+025}
\end{figure}

\begin{figure}[htb]
 \begin{center}
  {\includegraphics{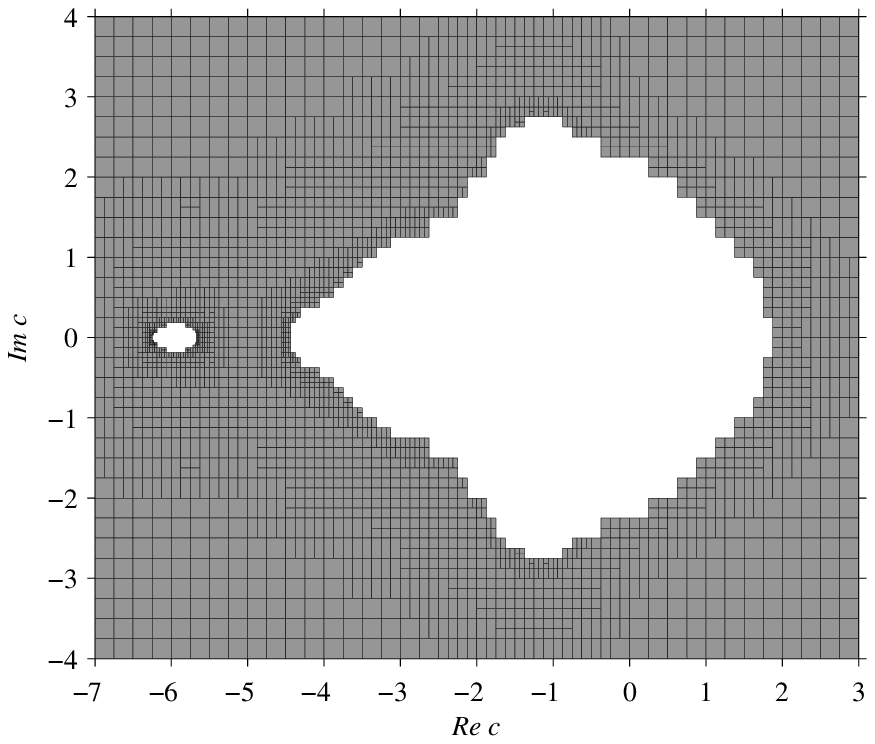}}
 \end{center}
 \caption{The shaded region is contained in 
 $\mathcal{H}^{\mathbb{C}} \cap \{a = -1\}$.}
 \label{FIG:ar} 
\end{figure}

\begin{figure}[htb]
 \begin{center}
  {\includegraphics{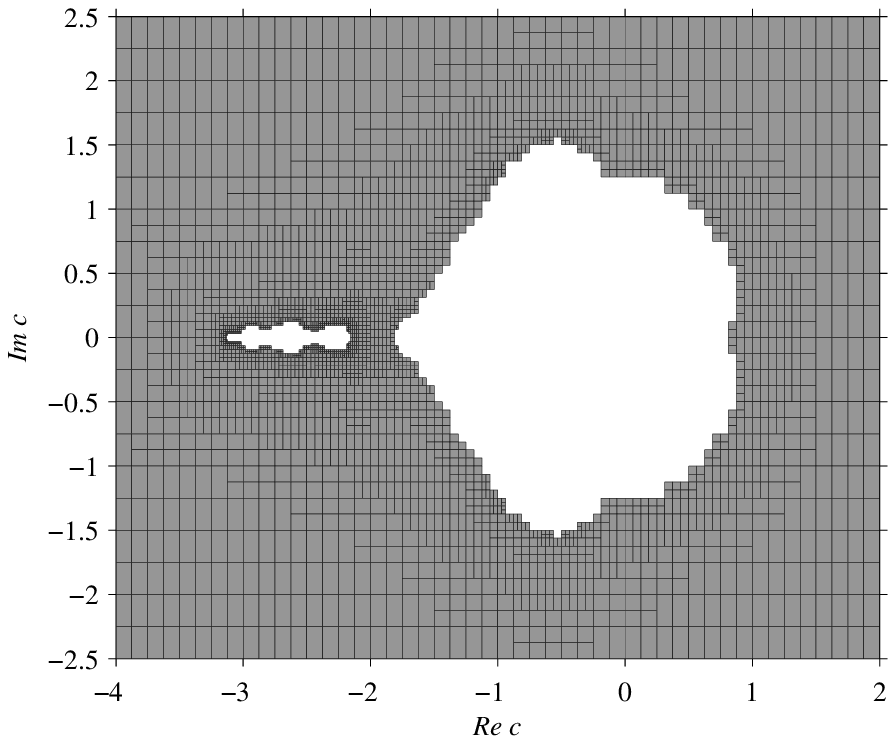}}
 \end{center}
 \caption{The shaded region is contained in 
 $\mathcal{H}^{\mathbb{C}} \cap \{a = -0.375\}$.}
 \label{FIG:a-0375}
\end{figure}

Now we discuss Theorem \ref{THM:existence}.
We begin by defining $I_p$, $I_q$, $I_r$ and $I_s$.
Let
\begin{align*}
 L_p &:= \{1\} \times \{\mathbb{C} \setminus 
 \text{white regions in Figure~\ref{FIG:ap}}\},\\
 L_q &:= \{0.25\} \times \{\mathbb{C} \setminus 
 \text{white regions in Figure~\ref{FIG:a+025}}\},\\
 L_r &:= \{-1\} \times \{\mathbb{C} \setminus 
 \text{white regions in Figure~\ref{FIG:ar}}\},\\
 L_s &:= \{-0.375\} \times \{\mathbb{C} \setminus 
 \text{white regions in Figure~\ref{FIG:a-0375}}\}
\end{align*}
and $L: = L_p \cup L_q \cup L_r \cup L_s$.
To be precise, 
these regions are defined by a finite 
number of closed rectangles.
The complete list of these rectangles is available at the author's web page.
The set $L_p \cap \mathbb{R}^2$ have three components:
two unbounded intervals, and one bounded interval.
We define $I_p$ to be this bounded one.
Similarly, $I_q$, $I_r$ and $I_s$ are defined to be the bounded intervals 
contained in $L_q \cap \mathbb{R}^2$, $L_q \cap \mathbb{R}^2$ and 
$L_s \cap \mathbb{R}^2$, respectively.

\begin{lemma}
 \label{LEMMA:hyp}
 If $(a, c) \in L$ then 
 $H_{a,c}$ is uniformly hyperbolic on its chain recurrent set 
 $\mathcal{R}(H_{a,c})$.
\end{lemma}
The proof of this lemma is computer assisted.
We leave it to \S \ref{SEC:hyperbolicity}.

Recall that the hyperbolicity of the chain recurrent set implies the
$\mathcal{R}$-structural stability \cite[Corollary 8.24]{Shub}.
Therefore, it follows from Lemma~\ref{LEMMA:hyp} that no bifurcation occurs
in $\mathcal{R}(H_{a,c})$ as long as $(a, c) \in L$.
Thus $\mathcal{R}(H_{a,c})$ is a hyperbolic full horseshoe for all 
$(a,c) \in L$.

Lemma~\ref{LEMMA:hyp} is not sufficient to conclude the hyperbolicity
of $K_{a,c}^{\mathbb{C}}$ because $\mathcal{R}(H_{a,c})$ and 
$K_{a,c}^{\mathbb{C}}$ do not necessarily coincide.
However, we can show that these sets are equal in the 
horseshoe locus, as follows.

\begin{corollary}
 \label{COR:hyp}
 If $(a, c) \in L$ then 
 $H_{a,c}|K_{a,c}^{\mathbb{C}}$ is a hyperbolic full horseshoe, that is, 
 $L \subset \mathcal{H}^{\mathbb{C}}$.
\end{corollary}

\begin{proof}[Proof of Corollary \ref{COR:hyp}]
 Let
 \begin{align*}
  K_{a,c}^+ :=& 
  \{p \in \mathbb{C}^2: \{H_{a,c}^n(p)\}_{n \geq 0} \text{ is bounded}\},\\
  K_{a,c}^- :=& 
  \{p \in \mathbb{C}^2: \{H_{a,c}^n(p)\}_{n \leq 0} \text{ is bounded}\}
 \end{align*} 
 and $J_{a, c}^{\pm} := \partial K_{a,c}^{\pm}$.
 Define $J_{a,  c} = J_{a, c}^+ \cap J_{a, c}^-$.
 Then $K_{a,c}^{\mathbb{C}} = K_{a, c}^+ \cap K_{a,c}^-$ and we have
 $J_{a,c} \subset \mathcal{R}(H_{a,c}) \subset K_{a, c}^{\mathbb{C}}$
 \cite[Proposition 9.2.6, Theorem 9.2.7]{MNTU}.
 Suppose $(a, c) \in L$.
 Since $\mathcal{R}(H_{a,c})$ is a full horseshoe, 
 all periodic points of $H_{a,c}$ is contained in $\mathcal{R}(H_{a,c})$
 and therefore they are of saddle type.
 Thus there exists no attracting periodic orbit.
 Furthermore, $J_{a,c}$ is uniformly hyperbolic because it is
 a closed sub-invariant set of $\mathcal{R}(H_{a,c})$.
 It follows that $\int K^+ = \emptyset$ \cite[Theorem~5.9]{BS91}.
 Since $\vert a \vert \leq 1$, we also have $\int K^- = \emptyset$
 \cite[Lemma 5.5]{BS91}.
 As a consequence, $J_{a,c}^+ = K_{a,c}^+$ and $J_{a,c}^- = K_{a,c}^-$,
 and hence $J_{a,c} = \mathcal{R}(H_{a,c}) = K_{a,c}^{\mathbb{C}}$.
 Therefore, Lemma~\ref{LEMMA:hyp} implies this corollary.
\end{proof}

From Corollary~\ref{COR:hyp} it follows that $I_p$, $I_q$, $I_s$ and $I_r$
are contained in $\mathcal{H}_0^{\mathbb{C}} \cap \mathbb{R}^2$.
To complete the proof of Theorem~\ref{THM:existence}, 
we need to show that these intervals are of type-3.

A simple and direct way for proving this is to
show that the number of periodic points contained in 
$K_{a, c}^{\mathbb{R}}$ is
non-zero and different from that of a full horseshoe.
Rigorous interval arithmetic and the Conley index theory can be applied 
for this purpose. 
We discuss this method in the appendix.

\begin{figure}[htb]
 \begin{center}
  \includegraphics{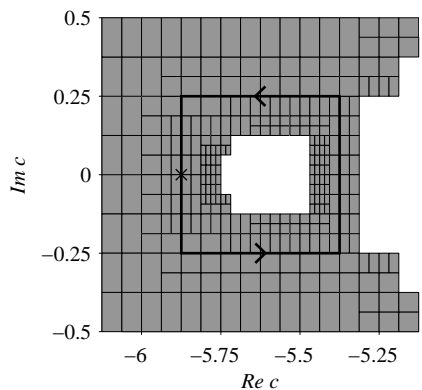}
 \end{center}
 \caption{The loop $\beta_p: [0,1] \to L_p$ based at
 $(a, c) = (1, -5.875)$.}
 \label{FIG:beta_p}
\end{figure}

Another way is to  make use of Theorem~\ref{THM:pruning}.
Since we have already shown that $L \subset \mathcal{H}_0^{\mathbb{C}}$,
we can consider the monodromy of loops in $L$, from which we
derive the information of $K_{a, c}^{\mathbb{R}}$.

Let $\beta_{p}: [0, 1] \to L_p$ be a loop that turns around
the smaller white island of Figure~\ref{FIG:ap}
as illustrated in Figure~\ref{FIG:beta_p}.
We require that $\beta_p(1/2) \in I_p$, and that $\beta_p$ be symmetric
with respect to the complex conjugation, that is, 
$\bar{\beta}_p = \beta_p^{-1}$.
Then we define a loop $\gamma_p: [0, 1] \to L_p \cup \mathrm{HOV}$
based at $(1, -10) \in \mathrm{DN}$ by setting
\begin{equation*}
 \gamma_p := \bar{\alpha}^{-1} \cdot \beta_p \cdot \alpha
\end{equation*}
where $\alpha:[0,1] \to \mathrm{HOV} \cup L_p$ is a path from 
$(1, -10)$ to the basepoint of $\beta_p$.
Choose the parametrization of $\gamma_p$ so that $\gamma_p(1/2) \in I_p$ and
$\bar{\gamma}_p = \gamma_p^{-1}$ hold. 
Similarly we define loops $\gamma_q$, $\gamma_r$ and $\gamma_s$
based at $(1, -10)$ turning around the smaller islands in $L_q$, $L_r$ 
and $L_s$, respectively.

\begin{proposition}
 \label{PROP:monodromy}
 The automorphism $\rho(\gamma_p)$ interchanges the words 
 $0010100$ and $0011100$ contained in 
 $s = (s_i)_{i \in \mathbb{Z}} \in \Sigma_2$.
 Namely, 
 \begin{equation*}
  (\rho(\gamma_p) (s))_i = \begin{cases}
                            0 & \text{if }
                            s_{i-3}\cdots s_i \cdots s_{i+3} = 0011100\\
                            1 & \text{if } 
                            s_{i-3}\cdots s_i \cdots s_{i+3} = 0010100\\
                            s_i & \text{otherwise.}
                           \end{cases}
 \end{equation*}
 Similarly, $\rho(\gamma_q)$ interchanges $10100$ and $11100$,
 $\rho(\gamma_r)$ interchanges $10010$ and $10110$, and
 $\rho(\gamma_s)$ interchanges $0010$ and $0110$. 
\end{proposition}

The proof of Proposition \ref{PROP:monodromy} is also computer assisted.
An algorithm for this will be discussed in \S \ref{SEC:monodromy}.

Now we are prepared to prove Theorem \ref{THM:existence}.
\begin{proof}[Proof of Theorem \ref{THM:existence}]
 Since $\Fix(\rho(\gamma_p))$ is a non-empty proper subset of $\Sigma_2$,
 Theorem~\ref{THM:pruning} implies that $\gamma_p(1/2) \in I_p$ is of type-3.
 By considering loops homotopic to $\gamma_p$, 
 we can show that all $(a,c) \in I_p$ are also of type-3.
 Proofs for other intervals are the same.
\end{proof}

Theorem \ref{THM:infinite_order} immediately follows from the following
proposition.

\begin{proposition}
 The order of $\psi = \rho(\gamma_{\emptyset}) \cdot \rho(\gamma_s)$
 is infinite.
\end{proposition}

The proof below is due to G. A. Hedlund~\cite[Theorem 20.1]{Hedlund}.

\begin{proof}
 For non-negative integer $p$, we define elements of $\Sigma_2$
 named $x^{(2p)}$ and $x^{(2p+1)}$ by
 \begin{align*}
  x^{(2p)} = \cdots 01010101{\mathbf{0110110}}(10)^{p}&.11111\cdots, \\
  x^{(2p+1)} = \cdots 01010101{\mathbf{0110110}}(10)^{p}1&.00000\cdots.
 \end{align*}
 We then look at the orbit of $x = x^{(0)}$ under the map $\psi$.
 A simple calculation shows that
 \begin{align*}
  x = \, \cdots01010101{\mathbf{0110110}}.11111\cdots& = x^{(0)},\\
  \psi(x) = \,  \cdots1010101{\mathbf{0110110}}1.00000\cdots& = x^{(1)}, \\
  \psi^2(x) = \,  \cdots010101{\mathbf{0110110}}10.11111\cdots& = x^{(2)}, \\
  \psi^3(x) = \,  \cdots10101{\mathbf{0110110}}101.00000\cdots& = x^{(3)}.
 \end{align*}
 By induction, it follows that $\psi^n(x^{(0)}) = x^{(n)}$.
 Since $x^{(n)} \ne x^{(m)}$ if $n \ne m$, this implies that
 the order of $\psi$ is infinite.
\end{proof}

Theorem~\ref{THM:SFT} is a direct consequence of 
Theorem~\ref{THM:pruning} and Proposition~\ref{PROP:monodromy}.

\section{Hyperbolicity}\label{SEC:hyperbolicity}
We recall an algorithm for proving
the uniform hyperbolicity of chain recurrent sets developed
by the author \cite{Arai05}.
We also refer the reader to the work of Suzanne Lynch Hruska 
\cite{Hru1, Hru2} for another algorithm.

Let $f$ be a diffeomorphism on a manifold $M$ and
$\Lambda$ a compact invariant set of $f$.
We denote by $T\Lambda$ the restriction of the tangent bundle
$TM$ to $\Lambda$.

\begin{definition}
 We say that $f$ is {\it uniformly hyperbolic} on $\Lambda$, or
 $\Lambda$ is a {\it uniformly hyperbolic invariant set} if 
 $T\Lambda$ splits into a direct sum $T\Lambda = E^s \oplus E^u$
 of two $Tf$-invariant subbundles and there exist constants $c > 0$
 and $0 < \lambda < 1$ such that $\|Tf^n|_{E^s}\| < c\lambda^n$ and 
 $\|Tf^{-n}|_{E^u}\| < c\lambda^n$ hold for all $n \geq 0$.
 Here $\|\cdot\|$ denotes a metric on $M$.
\end{definition}

Proving the uniform hyperbolicity of $f$ according to this usual 
definition is, in general, quite difficult.
Because we must control two parameters $c$ and $\lambda$ 
at the same time, and further,
we also need to constant a metric on $M$
adapted to the hyperbolic splitting.

To avoid this difficulty, we introduce a weaker notion of hyperbolicity
called ``quasi-hyperbolicity''.
We consider $Tf: T\Lambda \to T\Lambda$, 
the restriction of $Tf$ to $T\Lambda$, as a dynamical system.
An orbit of $Tf$ is said to be trivial if it is contained in
the image of the zero section.

\begin{definition}
We say that $f$ is {\it quasi-hyperbolic} on $\Lambda$ if 
$Tf:T\Lambda \to T\Lambda$ has no non-trivial bounded orbit.
\end{definition}

It is easy to see that hyperbolicity implies quasi-hyperbolicity.
The converse is not true in general.
However, when $f|_{\Lambda}$ is chain recurrent, these two notions 
are equivalent.

\begin{theorem}
 [\cite{CFS77, SS}]
 \label{THM:CFS}
 Assume that $f|_{\Lambda}$ is chain recurrent, that is, 
$\mathcal{R}(f|_{\Lambda}) = \Lambda$.
 Then $f$ is uniformly hyperbolic on $\Lambda$ if and only if 
 $f$ is quasi-hyperbolic on it.
\end{theorem}

The definition of quasi-hyperbolicity can be rephrased in terms of
isolating neighborhoods as follows.
Recall that a compact set $N$ is an isolating neighborhood 
with respect to $f$ if the maximal invariant set
\begin{equation*}
 \Inv(N, f) := \{x \in N \mid f^n(x) \in N \text{ for all } n \in \mathbb{Z}\}
\end{equation*}
is contained in $\int N$, the interior of $N$.
An invariant set $S$ of $f$ is said to be isolated if there is 
an isolating neighborhood $N$ such that $\Inv(N, f) = S$.

Note that the linearity of $Tf$ in fibers of $TM$ implies that 
if there exists a non-trivial bounded orbit of $Tf: T\Lambda \to T\Lambda$, 
then any neighborhood of the image of the zero-section
must contain a non-trivial bounded orbit.
Therefore, the definition of quasi-hyperbolicity 
is equivalent to saying that the image of the zero section of $T\Lambda$ 
is an isolated invariant set with respect to 
$Tf: T\Lambda \to T\Lambda$.
To confirm that $\Lambda$ is quasi-hyperbolic,
in fact, it suffice to find an isolating neighborhood
containing the image of the zero section.

\begin{proposition}[\cite{Arai05}, Proposition 2.5]
 \label{PROP:quasi-hyp}
 Assume that $N \subset T\Lambda$ is an isolating neighborhood
 with respect to $Tf: T\Lambda \to T\Lambda$ and
 $N$ contains the image of the zero-section of $T\Lambda$.
 Then $\Lambda$ is quasi-hyperbolic.
\end{proposition}

Next, we check that the hypothesis of Theorem~\ref{THM:CFS} is satisfied
in the case of the complex H\'{e}non map.
Let us define
\begin{gather*}
R(a,c) := \frac{1}{2}(1 + |a| + \sqrt{(1+|a|)^2 + 4c}),\\
S(a,c) := \{(x,y) \in \mathbb{C}^2: |x| \leq R(a,c), |y| \leq R(a,c)\}. 
\end{gather*}
Then the following holds as 
in the case of the real H\'{e}non map \cite[Lemma~4.1]{Arai05}.
\begin{lemma}
 \label{LEM:DN}
The chain recurrent set $\mathcal{R}(H_{a,c})$ is contained in $S(a,c)$.
Furthermore, $H_{a,c}$ restricted to $\mathcal{R}(H_{a,c})$ is chain recurrent.
\end{lemma}

To prove Lemma~\ref{LEMMA:hyp}, therefore, it suffice
to show that $\mathcal{R}(H_{a, c})$ is quasi-hyperbolic for $(a, c) \in L$.
By Proposition~\ref{PROP:quasi-hyp}, all we have to do
is to find an isolating neighbourhood that contains 
the image of the zero-section of $T\mathcal{R}(H_{a, c})$.
More precisely, it is enough to find $N \subset TM$ such that 
\[
 \mathcal{R}(H_{a,c}) \subset N \quad \text{and} \quad
 \Inv(N, TH_{a,c}) \subset \int N
\]
hold.
Here we identify $\mathcal{R}(H_{a,c})$ and its image by the 
zero-section of $TM$.
Since there are algorithms \cite[Proposition~3.3]{Arai05}
that efficiently compute rigorous outer approximations
of $\mathcal{R}(H_{a,c})$ and $\Inv (N, TH_{a,c})$,
these conditions can be checked on computers.

Now we fix the parameter $a$ to $+1$ (or $0,25$, $-0.375$, $-1$) 
and regard $\{H_{1,c}\}$ as a parametrized family 
with the parameter $c \in \mathbb{C}$.
In the parameter plane, we define 
\[
C := \{c \in \mathbb{C} : |\Im c | \leq 8 \text{ and } |\Re c| \leq 8 \}.
\]
If $c \not \in C$ then $(1, c) \in \mathrm{HOV}$, and thus
we do not need to check the hyperbolicity for such $c$.
Furthermore, our computation can be restricted to 
the case when $\Im c \geq 0$ because $H_{1, c}$ and $H_{1, \bar{c}}$ are
conjugate via $\phi$ and hence the hyperbolicity of these two
maps are equivalent.

Finally, we perform Algorithm~3.6 of \cite{Arai05}
for the family $\{H_{1, c}\}$ with the initial parameter set 
$C \cap \{\Im c \geq 0\}$.
The algorithm inductively subdivide the initial parameter set
and outputs a list of parameter cubes on which the quasi-hyperbolicity
is verified.
This proves the quasi-hyperbolicity of $\mathcal{R}(H_{1, c})$ for 
$(1, c) \in L_p$.
The quasi-hyperbolicity for $L_q$, $L_r$ and $L_s$ is
also obtained by applications of the same algorithm.

Performed on a 2.5GHz PowerPC G5 CPU, 
the computation takes
$1496.5$ hours, $1348.6$ hours, $1496.1$ hours and $1288.7$ hours
for $L_p$, $L_q$, $L_r$ and $L_s$,
respectively.

\section{Monodromy}\label{SEC:monodromy}
In this section, we develop an algorithm for computing the monodromy
homomorphism $\rho$.

Let $\gamma: [0, 1] \to \mathcal{H}^{\mathbb{C}}_0$ be a loop
based at $\gamma(0) = \gamma(1) = (a_0, c_0) \in \mathrm{DN}$.
Since $\rho(\gamma)$ is defined in terms of conjugacies
$h_t = h_{t}^{\gamma}: K_{\gamma(t)}^{\mathbb{C}} \to \Sigma_2$
along $\gamma$, we first discuss how to compute them.

\begin{figure}[ht]
 \begin{center}
  {\includegraphics{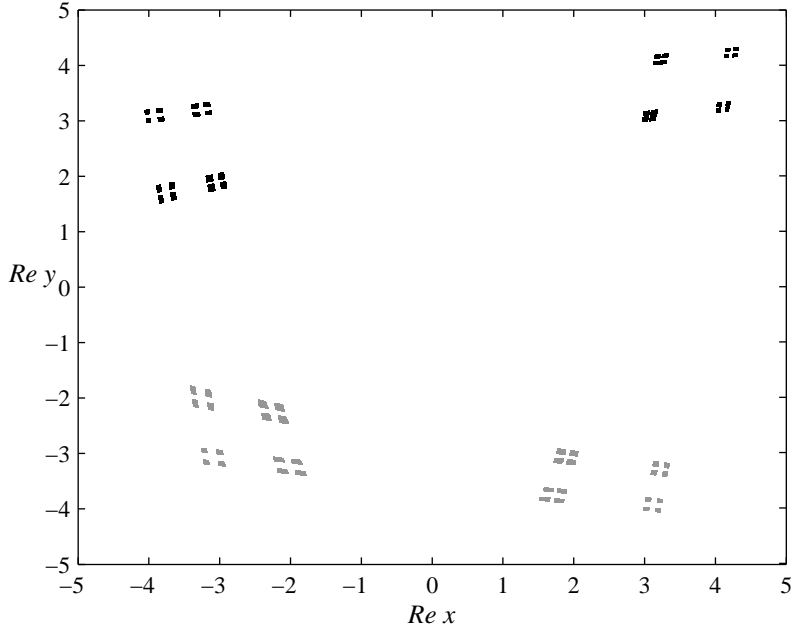}}
 \end{center}
 \caption{At $t = 0$: the initial partition $N_0^0$ and $N_0^1$.}
 \label{FIG:partition_0}
\end{figure}

\begin{figure}[ht]
 \begin{center}
  {\includegraphics{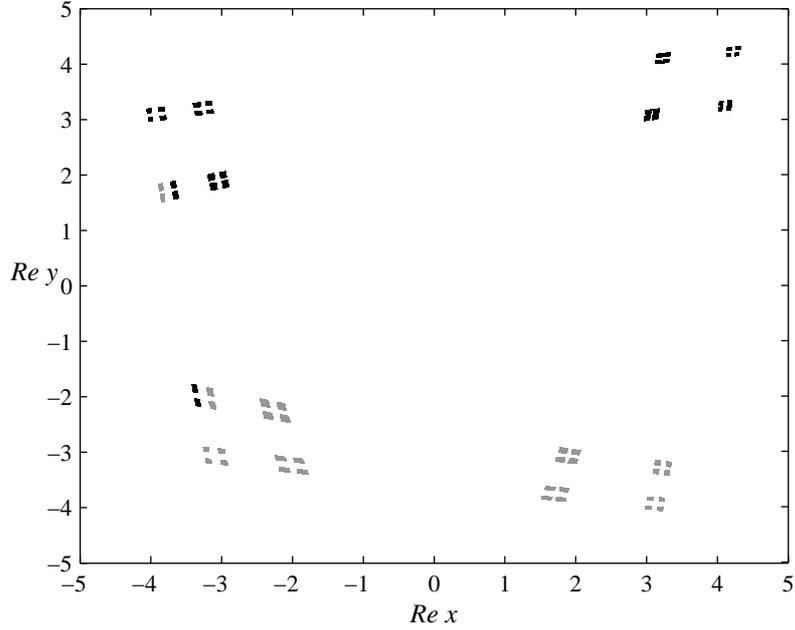}}
 \end{center}
 \caption{At $t = 1$: the partition $N_1^0$ and $N_1^1$,
 obtained by continuing $N_0^0$ and $N_0^1$along $\gamma_q$.}
 \label{FIG:partition_1}
\end{figure}

Let us recall the definition of $h_{t}$.
Define
\begin{equation*}
 K_0^0 := \{(x, y) \in K_{\gamma(0)}^{\mathbb{C}}: \Re y \leq 0\},\quad
 K_0^1 := \{(x, y) \in K_{\gamma(0)}^{\mathbb{C}}: \Re y \geq 0\}.
\end{equation*}
By the argument of Devaney and Nitecki \cite{DN79}, 
we have $K_0^0 \cap K_0^1 = \emptyset$ and the partition
$K_{\gamma(0)}^{\mathbb{R}} = K_0^0 \cup K_0^1$ 
induces a topological conjugacy $h_0$.
The continuation of this partition along $\gamma$ is defined by
\begin{gather*}
 K_t^0 := \{ z \in K_{\gamma(t)}^{\mathbb{C}}: 
  \text{ the continuation of } z 
  \text{ along } \gamma \text{ at } t = 0
  \text{ is in } K_0^0\}, \\
 K_t^1 := \{ z \in K_{\gamma(t)}^{\mathbb{C}}: 
  \text{ the continuation of } z 
  \text{ along } \gamma \text{ at } t = 0
  \text{ is in } K_0^1\}.
\end{gather*}
The conjugacy $h_t$ is, by definition, 
the symbolic coding with respect to this partition.
Namely, 
\begin{equation*}
 (h_t(z))_i := \begin{cases}
                0  & \text{if} \quad H_{\gamma(t)}^{i} (z) \in K_t^0\\
                1  & \text{if} \quad H_{\gamma(t)}^{i} (z) \in K_t^1.
               \end{cases}
\end{equation*}
To determine this conjugacy, however,
we do not need to compute $K_t^0$ and $K_t^1$ exactly.
It suffice to have rigorous outer approximations of them.
That is, if $N_t^0$ and $N_t^1$ are disjoint subsets of $\mathbb{C}^2$
such that $K_t^0 \subset N_t^0$ and $K_t^1 \subset N_t^1$ hold 
for all $t \in [0, 1]$,
then $k_{t}: K_{\gamma(t)}^{\mathbb{C}} \to \Sigma_2$
defined by
\begin{equation*}
 (k_t(z))_i := \begin{cases}
                0  & \text{if} \quad H_{\gamma(t)}^{i} (z) \in N_t^0\\
                1  & \text{if} \quad H_{\gamma(t)}^{i} (z) \in N_t^1
               \end{cases}
\end{equation*}
is identical to $h_t$.

Here is an algorithm to construct such $N_t^0$ and $N_t^1$. 

\begin{enumerate}
 \item[step 1.] Subdivide the interval $[0, 1]$ into $n$ closed intervals 
       $I_1, I_2, \ldots, I_n$ of equal length.

 \item[step 2.]
       Using interval arithmetic,
       we compute a cubical set $\mathcal{N}_k$
       for each $1 \leq k \leq n$
       such that $K^{\mathbb{C}}_{a,c} \subset \mathcal{N}_k$
       rigorously holds for all $(a, c) \in \gamma(I_k)$.
       Define $N_t := \mathcal{N}_k$ for $t \in I_k$.
 \item[step 3.] 
       Consider the set
       \begin{equation*}
        N := \bigcup_{t \in [0, 1]} \{t\} \times N_t \subset
         [0, 1] \times \mathbb{C}^2.
       \end{equation*}       
       Let $N^0$ and $N^1$ be the unions of the components of $N$
       which intersect with $\{0\} \times \{\Re y \leq 0\}$ and
       $\{0\} \times \{\Re y \geq 0\}$, respectively.
       If $N^0 \cap N^1 = \emptyset$, define $N_t^0 = N_t \cap N^0$ and
       $N_t^1 = N_t \cap N^1$ then stop.
       If this is not the case, we refine the subdivision of 
       $[0, 1]$ and the grid size of $\mathbb{C}^2$, and then go 
       back to step 1.
\end{enumerate}

\begin{figure}[htb]
 \begin{center}
  {\includegraphics{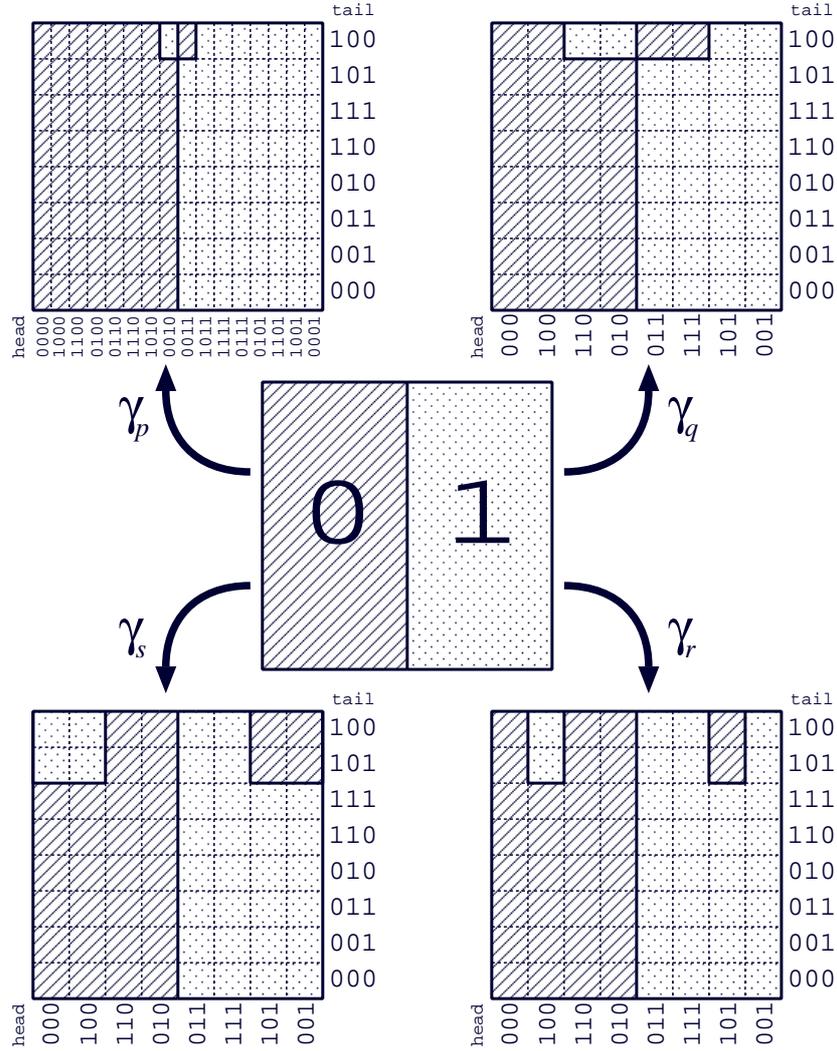}}
 \end{center}
 \caption{The change of the partition along 
 $\gamma_p, \gamma_q, \gamma_r$ and $\gamma_s$.}
 \label{FIG:monodromy}
\end{figure}

Applying the algorithm above to the loop $\gamma_q$, we obtain
Figure~\ref{FIG:partition_0} and \ref{FIG:partition_1}.
The interval $[0, 1]$ is decomposed into $n = 2^8$ sub-intervals, 
and the size of the grid for $\mathbb{C}^2$ is $2^{-8}$ in each 
direction.
The lightly and darkly shaded regions in Figure \ref{FIG:partition_0} 
are $N_0^0$ and $N_0^1$.
Similarly, Figure~\ref{FIG:partition_1} illustrates $N_1^0$ and $N_1^1$.
Notice that two partitions differ only in four blocks on 
the left hand side: two blocks of each of $N_0^0$ and $N_0^1$ are 
interchanged.
Using rigorous interval arithmetic, these blocks are identified as blocks
corresponding to the symbol sequences
$10.100$ and $11.100$ where the dot separates the
head and the tail of a sequence.
By the {\it head} of $s = (s_i)_{i \in \mathbb{Z}}$ we mean the sequence 
$\{\dots s_{-3} s_{-2} s_{-1} s_0\}$ and by the {\it tail}
$\{s_1 s_2 s_3 \dots\}$.

We execute the same computation also for loops $\gamma_p$, 
$\gamma_r$ and $\gamma_s$.
This yields Figure~\ref{FIG:monodromy}, which shows 
a schematic picture of the change along these loops.
Notice that ``head'' and ``tail'' labels in the figure 
indicates the symbol coding according to the initial partition
$N_0^0$ and $N_0^1$, illustrated in the central square.

Now, we can compute the image of $\gamma = \gamma_p$
(or $\gamma_q$, $\gamma_r$, $\gamma_s$) by $\rho$ as follows: 
Choose a symbol sequence $s = (s_i)_{i \in \mathbb{Z}} \in \Sigma_2$.
Then $z := h_0(s)$ is located in the central
square of Figure~\ref{FIG:monodromy}. 
By definition, $\rho(\gamma_p)(s)$ is the symbolic coding of
the same point $z$, 
but with respect to the partition on the top left corner of 
Figure~\ref{FIG:monodromy}.
Since two partitions differ only on blocks $0010.100$ and $0011.100$,
it follows that $(\rho(\gamma(s)))_i \ne s_i$ if and only if 
$H_{\gamma(0)}^i(z)$ is contained in these bocks.
Namely, 
\begin{equation*}
 (\rho(\gamma(s)))_i = 
  \begin{cases}
   0 &\text{ if } s_{i-3}s_{i-2}s_{i-1}s_{i}s_{i+1}s_{i+2}s_{i+3} = 0011100 \\
   1 &\text{ if } s_{i-3}s_{i-2}s_{i-1}s_{i}s_{i+1}s_{i+2}s_{i+3} = 0010100 \\
   s_i & otherwise.
  \end{cases}
\end{equation*}
Similarly we can compute $\rho(\gamma_q)$, 
$\rho(\gamma_r)$ and $\rho(\gamma_s)$.
This proves Proposition~\ref{PROP:monodromy}.

\appendix
\section*{Counting Periodic Orbits}

In this appendix, we prove Theorem~\ref{THM:existence} directly from 
Corollary~\ref{COR:hyp}, without any monodromy argument.
Instead of using Theorem~\ref{THM:pruning}, we show that the 
number of periodic points in $K_{a,c}^{\mathbb{R}}$ is different from 
that of a full horseshoe.
Specifically, we claim that the number of points in 
$\Fix(H_{a,c}^n) \cap \mathbb{R}^2$ is exactly as in Figure \ref{TABLE:pp}.

\begin{figure}[htb]
 \begin{center}
 \begin{tabular}{|c|c|c|c|c|c|c|c|c|} \hline
  & $\quad\mathrm{DN}\quad$ 
  & $L_p \cap \mathbb{R}^2$ 
  & $L_q \cap \mathbb{R}^2$ 
  & $L_r \cap \mathbb{R}^2$
  & $L_s \cap \mathbb{R}^2$
  & $\ \mathrm{EMP}\ $ \\ \hline 
  $n = 3$ &   8 &   8 &  8 &  2 &  2 & 0 \\
  $n = 4$ &  16 &  16 & 16 & 16 &  8 & 0 \\
  $n = 5$ &  32 &  22 & 22 & 22 & 12 & 0 \\
  $n = 6$ &  64 &  52 & 40 & 52 & 28 & 0 \\
  $n = 7$ & 128 & 114 & 72 & 72 & 44 & 0 \\ \hline
 \end{tabular}
 \end{center}
 \caption{The number of points in $\Fix(H_{a,c}^n) \cap \mathbb{R}^2$.}
 \label{TABLE:pp}
\end{figure}

We use the Conley index theory to prove the claim.
The reader not familiar with the Conley index may consult \cite{KMM,MM02}.

Assume $(a, c)$ is in one of $I_p$, $I_q$, $I_r$ or $I_s$.
We remark that 
the uniform hyperbolicity of $K_{a,c}^{\mathbb{R}}$ implies
that the number of real periodic points
is constant on these intervals.

First we compute a lower bound for the number of periodic points.
We begin with finding periodic points numerically.
Since periodic points are of saddle type and hence are numerically
unstable, we apply the subdivision algorithm \cite{GAIO} 
to find them.
For each periodic orbit found numerically,
we then construct a cubical index pair \cite{KMM}.
The existence of a periodic point in this index pair is then proved
by the following Conley index version of Lefschetz fixed point theorem.

\begin{theorem}[{\cite[Theorem~10.102]{KMM}}]
 \label{THM:Lefschetz}
 Let $(P_1, P_0)$ be an index pair for $f$ and $f_{P*}$ the homology index
 map induced by $f$.
 If $\sum_k (-1)^k \tr f^n_{P*k} \ne 0$
 then $\Inv(\cl(P_1 \setminus P_0),f)$ contains a fixed point of $f^n$.
\end{theorem}
This theorem assures that there exists at least one periodic orbit in each
index pair, and therefore we obtain a lower bound for the number of
points in $\Fix(H_{a,c}^n) \cap \mathbb{R}^2$.

To compute an upper bound, we have two methods.

One is to prove the uniqueness of the periodic orbit in each index pair.
As long as the size of the grid used in the subdivision algorithm was 
fine enough, we can expect that each index pair isolates exactly one
periodic orbit of period $n$.
Since periodic points are hyperbolic,
uniqueness can be achieved by a Hartman-Grobman type theorem 
\cite[Proposition~4.1]{AM}.

The other one is to use the fact that the number of fixed points of
$H_{a,c}^n: \mathbb{C}^2 \to \mathbb{C}^2$ is independent of the parameter,
in fact it is  $2^n$,
counted with multiplicity \cite[Theorem~3.1]{FM}.
In our case, 
the uniform hyperbolicity implies that the multiplicity is always $1$ and 
hence there are exactly $2^n$ distinct points in $\Fix(H_{a,c}^n)$.
Therefore, if we find $k$ distinct fixed points of $H_{a,c}^n$
outside $\mathbb{R}^2$, then the cardinality of
$\Fix(H_{a,c}^n) \cap \mathbb{R}^2$ 
must be less than or equal to $2^n -k$.
Again, we can apply Theorem~\ref{THM:Lefschetz} to establish the existence 
of fixed points in $\mathbb{C}^2 \setminus \mathbb{R}^2$.
This gives an upper bound.

For all cases shown in Figure~\ref{TABLE:pp},
the lower and upper bounds obtained by methods above coincide.
Thus our claim follows.


\begin{thebibliography}{99}
 \bibitem{Arai05}Z. Arai,
         {On Hyperbolic Plateaus of the H\'{e}non map}, 
         to appear in {\it Experimental Mathematics}.

 \bibitem{AM}Z. Arai and K. Mischaikow,
         {Rigorous computations of homoclinic tangencies},
         {\it SIAM Journal on Applied Dynamical Systems}
         {\bf 5} (2006), 
         280--292.

 \bibitem{BLS}E. Bedford, M. Lyubich and J. Smillie,
         {Polynomial diffeomorphisms of $C^2$. 
         IV: The measure of maximal entropy and laminar currents},
         {\it Invent. math.}
         {\bf 112} (1993),
         77--125.

 \bibitem{BS91}E. Bedford and J. Smillie,
         {Polynomial diffeomorphisms of $C^2$: 
         currents, equilibrium measure and hyperbolicity},
         {\it Invent. math.}
         {\bf 103} (1991),
         69--99.

 \bibitem{BS05}E. Bedford and J. Smillie,
         {The H\'{e}non family:
         The complex horseshoe locus and real parameter values},
         {\it Contemp. Math.}
         {\bf 396} (2006), 
         21--36.

 \bibitem{BDK1991}P. Blanchard, R. L. Devaney and L. Keen,
         {The dynamics of complex polynomials and automorphisms of the shift},
         {\it Invent. Math.}
         {\bf 104} (1991), 
         545--580.

 \bibitem{BLR1988}M. Boyle, D. Lind and D. Rudolph,
         {The automorphism group of a shift of finite type},
         {\t  Trans. Amer. Math. Soc.}
         {\bf 306} (1988),
         71--114.

 \bibitem{CFS77}R. C. Churchill, J. Franke and J. Selgrade,
         {A geometric criterion for hyperbolicity of flows},
         {\it Proc. Amer. Math. Soc.}
         {\bf 62} (1977), 
         137--143.

 \bibitem{Cvi91}P. Cvitanovi\'{c}, 
         {Periodic orbits as the skeleton of classical and quantum chaos},
         {\it Physica D}
         {\bf 51} (1991),  
         138--151. 

 \bibitem{DMS91}M. J. Davis, R. S. MacKay and A. Sannami,
         {Markov shifts in the H\'{e}non family},
         {\it Physica D}
         {\bf52} (1991), 
         171--178.

 \bibitem{GAIO}M. Dellnitz and O. Junge,
         {Set oriented numerical methods for dynamical systems},
         {\it Handbook of dynamical systems II},
         North-Holland, 2002, 221--264.

 \bibitem{DN79}R. Devaney and Z. Nitecki, 
         {Shift automorphisms in the H\'{e}non mapping},
         {\it Commun. Math. Phys.}
         {\bf 67} (1979),
         137--146.

 \bibitem{FM} S. Friedland and J. Milnor, 
         {Dynamical properties of plane polynomial automorphisms},
         {\it Ergodic Theory Dynam. Systems}
         {\bf 9} (1989), 67--99. 

 \bibitem{Hedlund}G. A. Hedlund,
         {Endomorphisms and automorphisms of the shift dynamical system},
         {\it Mathematical Systems Theory}
         {\bf 3} (1969), 320--375. 

 \bibitem{Hru1}S. L. Hruska,
         {A numerical method for constructing the hyperbolic structure 
         of complex H'{e}non mappings},
         {\it Found. Comput. Math.}
         {\bf 6} (2006),
         427--455.

 \bibitem{Hru2}S. L. Hruska,
         {Rigorous numerical models for the dynamics of
         complex H'{e}non mappings on their chain recurrent sets},
         {\it Discrete Contin. Dyn. Syst.}
         {\bf 15} (2006),
         529--558.

 \bibitem{KMM}T. Kaczynski, K. Mischaikow and M. Mrozek,
         {\it Computational Homology},
         Applied Mathematical Sciences, \textbf{157},
         Springer-Verlag, 2004.

 \bibitem{Kitchens}B. P. Kitchens,
         {\it Symbolic Dynamics},
         Springer-Verlag, 1998.
         
 \bibitem{MM02}K. Mischaikow and M. Mrozek,
         The Conley index theory, 
         {\it Handbook of Dynamical Systems II},
         North-Holland, 2002, 393--460.

 \bibitem{MNTU}S. Morosawa, Y. Nishimura, M. Taniguchi and  T. Ueda,
         {\it Holomorphic Dynamics},
         Cambridge University Press, 2000.

 \bibitem{SS}R. J. Sacker and G. R. Sell,
         {Existence of dichotomies and invariant splitting for 
         linear differential systems I},
         {\it J. Differential Equations}
         {\bf 27} (1974)
         429--458.

 \bibitem{Shub}M. Shub,
         {\it Global stability of dynamical systems},
         Springer-Verlarg, New York, 1987.
\end{thebibliography}
\end{document}